\documentclass[pdftex]{article}
\usepackage{amsmath,amssymb,latexsym,amsfonts}
\pdfoutput=1
\usepackage{amsthm}

\usepackage{caption}

\usepackage{color}
\usepackage{yhmath}
\usepackage{euscript}

\newcommand{\x}{\times}

\def\qed{\ensuremath{\qquad\hfill\Box}}

\newcommand{\R}{\ensuremath{\mathbb R}}
\newcommand{\C}{\ensuremath{\mathbb C}}
\newcommand{\N}{\ensuremath{\mathbb N}}

\newcommand{\m}{\operatorname{m}}

\newcommand{\al}{\ensuremath{\alpha}}
\newcommand{\be}{\ensuremath{\beta}}
\newcommand{\ga}{\ensuremath{\gamma}}
\newcommand{\de}{\ensuremath{\delta}}
\newcommand{\ep}{\ensuremath{\varepsilon}}

\newcommand{\la}{\ensuremath{\lambda}}

\newcommand{\ffi}{\ensuremath{\varphi}}

\newcommand{\om}{\ensuremath{\omega}}
\newcommand{\Om}{\ensuremath{\Omega}}

\newcommand{\La}{\ensuremath{\Lambda}}

\newcommand{\scrA}{\ensuremath{\EuScript{A}}}

\newcommand{\scrV}{\ensuremath{\EuScript{V}}}

\newcommand{\Cmn}{\ensuremath{{\mathbb C}^{m \times n}}}
\newcommand{\Cnn}{\ensuremath{{\mathbb C}^{n \times n}}}

\newcommand{\sumasec}[2]{\ensuremath{{#1}_1 + \cdots + {#1}_{#2}}}

\newcommand{\s}[1]{\ensuremath{\Lambda(#1)}}

\DeclareMathOperator{\rk}{rank}

\newcommand{\secsim}[2]{\ensuremath{#1_1,
\dots,#1_{#2}  }}  

\newcommand{\abso}[1]{\lvert #1 \rvert}

\newcommand{\norma}[1]{\lVert #1 \rVert}

\def\recoigu#1#2#3{\ensuremath{#1 = #2,\dots,#3}}

\def\demo{\noindent{\sc Proof.} }

\newcommand{\tresvari}[3]{#1_1(#2),\ldots,
#1_{#3}(#2)}

\def\recoper#1#2#3{\ensuremath{#1 \in \{#2,\dots,#3\} }}

\newcommand{\re}[2]{\ensuremath{#1,\ldots,#2}}

\newcommand{\secu}[2]{\ensuremath{#1_1,#1_2,\dots,
#1_{#2}}}

\newcommand{\reco}[3]{\ensuremath{#1 = #2,\dots,#3}}

\newcommand{\sumasecu}[2]{\ensuremath{#1_1 + #1_2 +\dots + #1_#2}}

\newcommand{\polimonivar}[4]{\ensuremath{#1^#3+#2_1(#4) #1^{#3-1}+\dots+#2_{#3-1}(#4)#1+ #2_#3(#4)}}

\newcommand{\card}{\#}

\newcommand{\secutres}[3]{#1_{#2_1},\ldots,#1_{#2_#3}}

\newcommand{\cuatrovari}[4]{#1_{#2_1}(#3),
\ldots,#1_{#2_#4}(#3)}

\newcommand{\tla}{\tilde{\lambda}}

\newcommand{\ttla}{\tilde{\tilde{\lambda}}}

\newcommand{\hla}{\hat{\lambda}}

\newcommand{\hhla}{\hat{\hat{\lambda}}}

\newtheorem{theorem}{Theorem}
\newtheorem{lemma}[theorem]{Lemma}
\newtheorem{proposition}[theorem]{Proposition}
\newtheorem{corollary}[theorem]{Corollary}

\newtheorem{remark}[theorem]{Remark}

\theoremstyle{definition}
\newtheorem{definition}[theorem]{Definition}

\usepackage[pagebackref,citecolor=blue,colorlinks,urlcolor=blue]{hyperref}

\title{Multiplicities of the structured pseudoeigenvalues}
\author{Juan-Miguel Gracia\thanks{%
The University of the Basque Country,
Department of Applied Mathematics and Statistics,
7, Paseo de la Universidad,
01006 Vitoria-Gasteiz, Spain.  \texttt{juanmiguel.gracia@ehu.es}}}
\date{August 20, 2009}

\begin{document}
\maketitle
\begin{abstract}

The structured pseudospectra of a matrix $A$ are sets of complex numbers
that are eigenvalues of matrices $X$ which are near to $A$ and  have the
same entries as $A$  at a fixed set of places.
	The sum of multiplicities of the eigenvalues of $X$ inside each connected component of the structured pseudospectra of $A$ does not depend on $X$. This fact is known, but not so much as it should be. For this reason, we give here an elementary and detailed proof of the result.
\end{abstract}


\section{Introduction}

In this paper we will use the Euclidean norm for vectors and the spectral norm for matrices. Let $M\in\Cnn$ and we denote  its spectrum by $\s M$; for any $\mu\in \s M$, we denote by $\m(\mu,M)$ the algebraic multiplicity of $\mu$ as an eigenvalue of $M$. The nullity  of a matrix $M\in\Cmn$
is defined as $\nu(M):=n-\rk M$.
 For all complex number $z_0$ and any positive real $r$, we denote by $B(z_0,r)$ the open ball centered at $z_0$ and radius $r$.

 Let $A\in\Cnn$ and $\s A=\{\secsim \la p\}$; let $\eta$ be a real number such that
\[
0<\eta <\min_{i\ne j} \frac{\abso{\la_i-\la_j}}{2}. 
\]
By the continuity of the eigenvalues there exists a number $\de>0$ such that for all $A'\in\Cnn$ that satisfies $\norma{A'-A}<\de$ we have: (1) $\s{A'}\subset \cup_{i=1}^p B(\la_i,\eta)$; (2) for every $\recoigu i1p,$ 
\begin{equation}\label{eq:local}
\sum_{\al \in \s{A'}\cap B(\la_i,\eta)}\m(\al,A')=\m(\la_i,A).
\end{equation}
 Let $\ep >0$, we define the $\ep$-\emph{pseudospectrum} of $A$ as the set
$
\La_\ep(A):=\{z\in\C: z\in\s X , X\in\Cnn, \norma{X-A}\le \ep\}.
$
 If we change $\le$ by $<$ in this definition we define the \textit{strict} $\ep$-\emph{pseudospectrum} of $A$, denoted by
$
\La'_\ep(A).
$
We will call $\ep$-pseudoeigenvalue of $A$ to any complex number $\la$ of the $\ep$-pseudospectrum of $A$. When we need to precise some matrix $X,\norma{X-A}\le \ep$, such that $\la\in\s X$, we will say that $\la$ is an $\ep$-pseudoeigenvalue of $A$ \emph{associated with} $X$.

The number of elements of a finite set $Z$ will be denoted by $\# Z$.
Let $S$ be a subset of \underbar{$n$}$\x $\underbar{$n$} where \underbar{$n$}:=$\{1,2,\ldots,n\}$. Let $s:=\# S$ and we  denote by $(i_1,j_1),(i_2,j_2),\ldots,(i_s,j_s)$ the elements of $S$ ordered according the lexicographical order of $\N\x\N$, with $\N:=\{1,2,\ldots\}$ the set of natural numbers. Let us define the matrix function
\[
M_S\colon \C^s\to\Cnn
\]
that to each $s$-tuple $z=(\secsim zs)\in\C^s$ associates the matrix $M_S(z)$, given by
\[
  \begin{cases}
  M_S(z)(i_k,j_k)=z_k    &   \text{for  $k=1,\ldots,s$,}\\
M_S(z)(i,j)=0	  &   \mbox{for $(i,j)\notin S.$}
    \end{cases}
\]
For every $\ep >0$ we will call \emph{structured $\ep$-pseudospectrum of (structure $S$) of the matrix $A$} to the set
\[
\La_{S,\ep}(A):=\bigcup_{\substack{z\in\C^s\\ \norma{z}\le \ep}}\La(A+M_S(z));
\]
The \emph{strict} structured $\ep$-pseudospectrum of (structure $S$) of the matrix $A$ will be the set
\[
\La'_{S,\ep}(A):=\bigcup_{\substack{z\in\C^s\\ \norma{z}< \ep}}\La(A+M_S(z)).
\]
We will call structured $\ep$-pseudoeigenvalue of (structure $S$) of $A$ to any complex number $\la$ of $\La_{S,\ep}(A).$
Analogously, when we need to precise some $z,\norma{z}\le \ep$, such that $\la\in\s{A+M_S(z)}$, we say that $\la$ is an structured $\ep$-pseudoeigenvalue of $A$ \emph{associated with } $z$.

Embree and Trefethen~\cite{embree_tr:2001} have extended a large number of the properties of the spectrum to the pseudospectrum.
The main objective of this paper is to give a generalization of  equality~\eqref{eq:local}. In fact, Mosier showed~\cite[Theorem 2, p. 269]{mosier1986root} that for sufficiently small $\ep$, if $\norma{A'-A}<\ep$, then both $A$ and $A'$ have the same number of eigenvalues in each connected component of $\La'_\ep(A)$ (counted with multiplicity), and also that there is at least one eigenvalue of $A$ in each connected component. In his brilliant thesis about structured pseudospectra~\cite{karow2003geometry}, Karow
 gave a rigorous   proof of this result; see his Proposition 2.2.3. A sketchy proof, for matrix polynomials, was given by Lancaster and Psarrakos~\cite[Theorem 2.3]{lancaster2005pseudospectra}.     

This paper is organized as follows: Section 1 is this Introduction.
Section 2 reviews known properties about the compactness and closure of pseudospectra. Section 3 gives other proof of the fact that every connected component of $\La_{S,\ep}(A)$ contains an eigenvalue of $A$. In Section 4 we present formulas that give the number of $k$-multiple roots of a polynomial in terms of its coefficients. Several results on the parameterization  of the spectrum of a continuous matrix function of a real variable, are reformulated in Section 5. The main result about the constancy of the sum of multiplicities of the eigenvalues of $A+M_S(z),\norma{z}<\ep,$ inside a connected component of the structured $\ep$-pseudospectrum, is given in Section 6. A lower bound of the distance from a matrix $A$ to the set of matrices $X$ that have eigenvalues with greater multiplicity, is given in Section 7, jointly with other corollaries.

\section{Compactness and closure} 
\smallskip

\begin{proposition}\label{pro:compact}
Let $\ep>0$, then,
\renewcommand{\theenumi}{\arabic
{enumi}}
\renewcommand\labelenumi{
\upshape\theenumi .}
\begin{enumerate}
\item $\La_{S,\ep}(A)$ is a compact set.
\item The set of isolated points of $\La_{S,\ep}(A)$ is contained in $\s A$.
\item $\La'_{S,\ep}(A)\setminus \s A$ is a nonempty open set.
\end{enumerate}
\end{proposition}

 We denote by $\overline{X}$ the closure of the subset $X$ of $\C$.
\begin{corollary}\label{cor:clausura}
\[
\overline{\La'_{S,\ep}(A)}=\La_{S,\ep}(A).
\]
\end{corollary}

 See~\cite[p. 17--18, Propositions 2.2.3, 2.2.4 and 2.2.5]{karow2003geometry} and~\cite{burke2003optimization}.

\section{Connected components}
\smallskip

\begin{theorem} \label{th:uno}
Let $\ep >0$, for each $z\in\C^s$ such that $\norma{z}\le\ep$,  every connected component $T$ of $\La_{S,\ep}(A)$ contains at least an eigenvalue of $A+M_S(z)$.
\end{theorem}

\demo Let $w_0$ be any point of $T$; so, there is a vector $z^0\in\C^s$ such that $\norma{z^0}\le\ep$ and $w_0$ is an eigenvalue of $A+M_S(z^0)$.
  Let $u(t):=z^0+t(z-z^0),\; t\in [0,1]$. As any ball in a normed  space is convex~\cite[p. 356, Exercise 2]{lancaster1985tismenetsky}, for all $t\in[0,1]$, $u(t)\in \overline{B(0,\ep)}$; therefore, 
     for all $t\in[0,1]$, 
\[
\s{A+M_S(u(t))}\subset \La_{S,\ep}(A).
\] 
By~\cite[p. 155, Corollary VI.1.6]{bhatia:1997},\cite[p. 126, Theorem 5.2]{kato:1982} there exist continuous functions
$
\secsim \la n\colon [0,1]\to\C
$
such that for every $t\in [0,1]$,
\[
\s{A+M_S(u(t))}=\{\tresvari \la tn\}\quad \text{(it may have repetitions)}.
\]
Then,
\[
\s {A+M_S(z^0)}=\{\tresvari \la 0n\};
\]
hence, there is an index
 $\recoper i1n$
such that $w_0=\la_i(0)$. As
\[
\s {A+M_S(z)}=\{\tresvari \la 1n\},
\]
 $\la_i(1)\in \s {A+M_S(z)}$, and given that there exists a continuous path
$\la_i(t)$ that connects 
$w_0=\la_i(0)$ with $\la_i(1)$ inside  $\La_{S,\ep}(A)$, 
we have
$\la_i(1)\in T$. Hence, $T$ contains an eigenvalue of $A+M_S(z)$.\qed

\begin{remark}\label{obs:particular}
{\em In particular, the theorem that we have just demonstrated implies that every connected component of the  pseudospectrum, $\La_{S,\ep}(A)$, contains at least an eigenvalue of $A$, because $M_S(0)=0$. Other authors prove first this assertion from the maximum modulus principle for analytic functions,~\cite{burke2003optimization},~\cite{embree_tr:2001}.   
}
\end{remark}

\section{Number of roots}

All information concerning a polynomial is enclosed in its coefficients; so it must be possible to answer many questions about the polynomial from conditions satisfied by its coefficients.
In this section we will give thus the number of distinct roots that has a polynomial, or the number of distinct roots that have a prescribed multiplicity. This will be done by computing the rank of matrices constructed with the coefficients of the polynomial.

Consider the complex polynomial
\[
a(\la)=\la^n+a_1\la^{n-1}+\cdots+a_{n-1}\la+a_n
\]
and $h$ complex polynomials $b_1(\la),\ldots,b_h(\la).$
 Let 
 \[
 \Pi:=\{a(\la),b_1(\la),\ldots,b_h(\la)\}.
 \]
  Let us suppose that the greatest degree of the polynomials
  $$b_1(\la),\ldots,b_h(\la)$$ is $p\le n$.
  
   Then we write
\begin{gather*}
b_1(\la)=b_{1,n-p}\la^p+b_{1,n-(p-1)}\la^{p-1}+
\cdots +b_{1,n-1}\la+b_{1n}\\
\vdots \qquad \qquad \qquad \qquad \qquad \qquad \vdots\\
b_i(\la)=b_{i,n-p}\la^p+b_{i,n-(p-1)}\la^{p-1}+
\cdots +b_{i,n-1}\la+b_{in}\\
\vdots \qquad \qquad \qquad \qquad \qquad \qquad \vdots\\
b_h(\la)=b_{h,n-p}\la^p+b_{h,n-(p-1)}\la^{p-1}+
\cdots +b_{h,n-1}\la+b_{hn}
\end{gather*}
where $b_{i,n-p}\neq 0$ for at least an
$i\in\{\re{1}{h}\}$. Let us define the matrix $p\times (n+p)$ associated with $a(\la )$:
\[
S_0:=\begin{bmatrix}
 1 & a_1  & a_2 & \cdots & a_n & 0 & \cdots & 0 \\
0 & 1 & a_1  & \cdots& a_{n-1} & a_n  & \cdots & 0\\
\vdots  & & \ddots & \ddots & & \ddots & \ddots & \vdots \\  
0 & 0 & \cdots & 1 & a_1 & \cdots & a_{n-1} & a_n
   \end{bmatrix}
\]
and a matrix $n\times(n+p)$ associated with each $b_i(\la)$
\[
S_i:=\begin{bmatrix}
 b_{i,n-p} & b_{i,n-(p-1)}  &  \cdots & b_{in} & 0 & \cdots & 0 \\
0 & b_{i,n-p} & \cdots& b_{i,n-1} & b_{in}  & \cdots & 0 \\
\vdots  &   & \ddots & & \ddots & \ddots &  \\  
0 & 0 & \cdots & b_{i,n-p} & \cdots & b_{i,n-1} & b_{in}
   \end{bmatrix}
\]
for $i=\re{1}{h}$. The generalized Sylvester's resultant of the polynomials $$a(\la ),b_1(\la),\ldots,b_h(\la )$$ is defined in this way

\begin{equation}\label{eq:resultante_gen}
R:=R(\Pi):=R\big(a,b_1,\ldots,b_h):=\begin{bmatrix}
 S_0   \\
 S_1\\
 \vdots\\
 S_h
   \end{bmatrix}\in\C^{(nh+p)\times (n+p)}.
\end{equation}
stacking the matrices $S_0,S_1,\ldots,S_h$. 
With these notations we have the next theorem, proved in~\cite{barnett1983polynomials}. 

\smallskip

\begin{theorem} \label{teo:barnett}
The number of common roots of the polynomials
\[
a(\la ),b_1(\la),\ldots,b_h(\la),
\]
counting multiplicities in their greatest common divisor, is equal
to the nullity of the matrix $R$.
\end{theorem}

\smallskip

Given the complex polynomial
\[
f(\la)=\la^n+a_1\la^{n-1}+\cdots+a_{n-1}\la+a_n,
\]
for each integer $k$, $1\le k\le n$, let us denote by
$N_k$ (resp. $\rho _k$) the number of distinct roots of
 $f(\la)$ of multiplicity $\ge k $ (resp. $=k$).

\smallskip

\begin{theorem} 
Let $R(f,f',\ldots,f^{(k)})$ be the generalized Sylvester's resultant matrix of the polynomials $f(\la), f'(\la),\ldots,f^{(k)}(\la)$, for each $k=1,2,\ldots,n,n+1$.
	Then,
	\begin{gather*}
N_k=\nu \big(R(f,f',\ldots,f^{(k-1)})\big)-\nu \big(R(f,f',\ldots,f^{(k)})\big)\\
\rho_k=\nu \big(R(f,f',\ldots,f^{(k-1)})\big)-2\nu \big(R(f,f',\ldots,f^{(k)})\big)+\\
\nu \big(R(f,f',\ldots,f^{(k+1)})\big)
\end{gather*}
for $1\le k\le n$, where if $k=1$,
$\nu \big(R(f,f',\ldots,f^{(k-1)})\big):=\nu \big(R(f)\big)$ and
\[
R(f):=\begin{bmatrix}
 1 & a_1  & a_2 & \cdots & a_n & 0 &\cdots & 0\\
 0 & 1 & a_1 & \cdots & a_{n-1} & a_n & \cdots & 0\\
 \vdots & & \ddots & \ddots & &&\ddots &\vdots\\
 0 & 0 & \cdots & 1 & a_1&  \cdots  &a_{n-1} & a_n
   \end{bmatrix}.
\]
\end{theorem}

\demo 
The matrix $R(f)$ has dimensions $(n-1)\times(2n-1)$;hence, $\rk R(f)=n-1$. Which implies 
\begin{equation}\label{eq:25-20-5}
\nu \big(R(f)\big)=n.
\end{equation} 

We will prove the formula for $N_k$ by induction on $k$.
Let us suppose that $k=1$. Then the number of distinct roots, $r$, of $f(\la)$ coincides with $N_1$. 
 Let $ \secu \la r \in \C$ be the distinct roots of $f(\la)$ with
 respective  multiplicities  $\secu mr$. Then we deduce that the number of common roots to $f(\la)$ and $f'(\la)$ counting their multiplicities in the $\gcd (f(\la),f'(\la))$, is
\begin{equation}\label{eq:25-20-20}
(m_1-1)+(m_2-1)+\cdots +(m_r-1),
\end{equation}
because the $m_i$-multiple root $\la_i$ of $f(\la)$ is an $(m_i-1)$-multiple root of $f'(\la )$, and, so, is a root of multiplicity  $m_i-1$ of $\gcd(f(\la),f'(\la ))$, for each $\reco i1r$. Therefore,
\[
(m_1-1)+(m_2-1)+\cdots +(m_r-1)=\nu \big(R(f,f')\big).
\]
Hence,
\[
\sumasecu mr -r=\nu \big(R(f,f')\big),
\]
but $\sumasecu mr=n$; which implies $n-r=\nu \big(R(f,f')\big)$; 
that is to say, 
 $r=n-\nu \big(R(f,f')\big)$. Therefore,
\[
r=\nu \big(R(f)\big)-\nu \big(R(f,f')\big);
\]
which is equivalent to 
\begin{equation}\label{eq:25-20-30}
N_1=\nu \big(R(f)\big)-\nu \big(R(f,f')\big),
\end{equation}
because $\nu \big(R(f)\big)=n$. This proves the formula for $N_k$ for $k=1$.

Next, let us suppose that the formula for $N_k$ is true for $1,2,\ldots,k$.
 Let $\secu{\al}{N_{k+1}}$ be the roots of $f(\la)$ with multiplicities  $\ge k+1$; let us denote their respective
multiplicities by 
$\secu{m}{N_{k+1}}$. Then for each 
$\recoigu i1{N_{k+1}}$, the number $\al_i$ is an $(m_1-(k+1))$-multiple root of the polynomial
\[
\gcd \big(f(\la),f'(\la),\ldots, f^{(k+1)}(\la)\big);
\]
this implies
\[
(m_1-(k+1))+\cdots +(m_{N_{k+1}}-(k+1))=\nu \big(R(f,f',\ldots ,f^{(k+1)})\big).
\]
Therefore,
\begin{equation}\label{eq:25-20-40}
\sumasec {m}{N_{k+1}}-(k+1)N_{k+1}=\nu \big(R(f,f',\ldots, f^{(k+1)})\big).
\end{equation}

But, by the definition of the numbers $N_j,\recoigu j1n$, the polynomial $f(\la )$ has $N_k-N_{k+1}$  $k$-multiple distinct roots, $N_{k-1}-N_k$ $(k-1)$-multiple distinct roots, \ldots, $N_2-N_3$ double distinct roots and $N_1-N_2$ simple distinct roots. Hence, 
\begin{align*}
 n&=(\sumasec m{N_{k+1}}) +k(N_k-N_{k+1})+(k-1)(N_{k-1}-N_k)\\
 & + \cdots + 2(N_2-N_3)+(N_1-N_2);      
\end{align*}
hence,
\begin{align}\label{eq:25-20-50}
 \sumasec m{N_{k+1}} &= n-kN_k+kN_{k+1}-(k-1)N_{k-1}\\ \nonumber
&\quad +(k-1)N_k-(k-2)N_{k-2}+(k-2)N_{k-1}\\ \nonumber
& \quad \cdots -2N_2+2N_3-N_1+N_2\\ \nonumber
&\!\!\!\! =  n+kN_{k+1}-N_k-N_{k-1}-\cdots -N_3-N_2-N_1;      
\end{align}
Now, \eqref{eq:25-20-40} and~\eqref{eq:25-20-50} imply
\begin{align*}
 n+kN_{k+1}-N_k-N_{k-1}-\cdots -N_2-N_1-(k+1)N_{k+1}\\
     = \nu \big(R(f,f',\ldots, f^{(k+1)})\big);
\end{align*}
therefore,
\[
n-\sum_{i=1}^k N_i - N_{k+1}=\nu \big(R(f,f',\ldots, f^{(k+1)})\big);
\]
or,
\begin{equation}\label{eq:25-20-60}
n-\sum_{i=1}^k N_i -\nu \big(R(f,f',\ldots, f^{(k+1)})\big) =N_{k+1}.
\end{equation}
But, by induction hypothesis,
\begin{align*}
 N_1 &= \nu \big(R(f)\big)-\nu \big(R(f,f')\big), \\
 N_2 &= \nu \big(R(f,f')\big)-\nu \big(R(f,f',f'')\big), \\
& \vdots \\
 N_{k-1} &= \nu \big(R(f,f',\ldots,f^{(k-2)})\big)-\nu \big(R(f,f',\ldots,f^{(k-1)})\big),  \\ 
N_k &= \nu \big(R(f,f',\ldots,f^{(k-1)})\big)-\nu \big(R(f,f',\ldots,f^{(k)})\big);   
\end{align*}
that implies
\[
\sum_{i=1}^kN_i=\nu\big(R(f)\big)-\nu \big(R(f,f',\ldots,f^{(k)})\big).   
\]
From here, by~\eqref{eq:25-20-60} and~\eqref{eq:25-20-5},
\[
\nu \big(R(f,f',\ldots,f^{(k)})\big)-\nu \big(R(f,f',\ldots,f^{(k+1)})\big)=N_{k+1}.
\]
In consequence, the formula for $N_k$ is true for all $\recoigu k1n$.

When $1\le k <n$, as $\rho_k:=N_k-N_{k+1}$,
\begin{multline*}
\rho_k=\nu \big(R(f,f',\ldots,f^{(k-1)})\big)-\nu \big(R(f,f',\ldots,f^{(k)})\big)\\
-\nu \big(R(f,f',\ldots,f^{(k)})\big)+\nu \big(R(f,f',\ldots,f^{(k+1)})\big).
\end{multline*}
If $k=n$, then $\nu \big(R(f,f',\ldots,f^{(n)})\big)=0$,
and $N_{n+1}=0$. Whence,
\[
\rho_n=N_n-N_{n+1}=N_n=\nu \big(R(f,f',\ldots,f^{(n-1)})\big)-\nu \big(R(f,f',\ldots,f^{(n)})\big);
\]    
or
\[
\rho_n=\nu \big(R(f,f',\ldots,f^{(n-1)})\big)-2\cdot 0+0.
\]
Therefore, the formula for $\rho_k$ is true for all $\recoigu k1n$. \qed

\smallskip

\begin{corollary}\label{cor:distintas}
The number of distinct roots of $f(\la)$ is
\begin{equation}\label{eq:distintas}
u=n-\nu\big(R(f,f')\big).
\end{equation}
\end{corollary}

\smallskip

The formula~\eqref{eq:distintas} appeared in the page 609 of \cite{gohberg_la_ro:1986}, just before (19.2.1).

\section{Spectrum of matrix functions of a real variable}

First, let us cite a lemma about locally constant functions on connected spaces.

\smallskip

\begin{lemma} \label{lem:local_constant}
Let $\Om$ be a connected topological space,  let $S$ be any set, and let $f:\Om\to S$ be a function such that for every $\om_0$, there is a neighborhood $\scrV$ of $\om_0$ such that the restriction of $f$ to $\scrV$ is constant. Then $f$ is constant on $\Om$.
\end{lemma}

\smallskip

For a proof see Lemma 2.1, p. 152 of~\cite{evard:1985}. Second, we write a lemma about the continuity of the spectrum of a continuous matrix function.

\smallskip
\begin{lemma} \label{lem:cont_espectro}
Let $E$ be a topological space and $\scrA  \colon E\to \Cnn $ a continuous matrix function. Let $t_0\in E, \secsim \al p$ be the distinct eigenvalues of $\scrA(t_0)$, and let $\eta$ be a real number such that
\[
0<\eta <\min_{i\ne j} \frac{\abso{\al_i-\al_j}}{2}. 
\]
Then there is a neighborhood $\scrV (t_0)$ of $t_0$ such that, for all $t\in\scrV(t_0)$,
\[
\La(\scrA(t))
\subset B(\al_1,\eta) \cup\cdots \cup B(\al_p,\eta),
\]
and for each $\recoper k1p,$
\[
\sum_{\mu\in\La(\scrA(t))\cap B(\al_k,\eta)}\m(\mu,\scrA(t))=\m(\al_k,\scrA(t_0)).
\]
\end{lemma}

\smallskip

\begin{definition}\label{bifur}
Let $\scrA:E\to \Cnn$ a continuous matrix function defined on a topological space $E$. Let
\[
N(t):=\# \{\mu\in\C: \det(\mu I-\scrA(t))=0\}
\]
for all $t\in E$. We will say that $t_0\in E$ \emph{is a bifurcation point of the spectrum of} $\scrA$ if the function $N$ is not constant on every neighborhood of $t_0$.
\end{definition}

\subsection*{Condition for a set of continuous functions not to have bifurcation points}

\smallskip

\begin{lemma} \label{lem:nobifurca}
Let $\Om$ be a connected topological space. Let $S$ be a set of continuous complex functions defined on $\Om$ such that the number of elements of the subset $S(t):=\{f(t)| f\in S\}$ of $\C$ is finite and does not depend on $t\in\Om$. 
  Let $t_0\in\Om$ and $f,g\in S$ be such that $f(t_0)=g(t_0)$. Then $f=g$.
\end{lemma}

See a proof in~\cite[Lemma 4.1, p. 375]{evard_gr:1990}.
As a consequence of Lemma~\ref{lem:nobifurca} we have the following lemma.

\smallskip

\begin{lemma} \label{lem:maximo}
Let $A\colon I \to \Cnn$  be a continuous matrix function defined on an interval $I$ of the real  line. Let us denote
\[
r:=\max_{t\in I}\#\, \s{A(t)}.
\] 
Then there are $r$ continuous functions
 $\la_j\colon I \to\C, \reco j1r$, such that for each $t\in I$,
\[
\{\tresvari \la t r\}=\s{A(t)}.  
\]
\end{lemma}

\smallskip

Let $A\colon (a,b)\to \Cnn$
be a continuous matrix function defined on an open interval $(a,b)$
of the real line. First, we are going to prove the claim:
\textit{An accumulation point of bifurcation points is also a bifurcation point. }

Let $\tau_0 \in (a,b)$ be an accumulation point of bifurcation points of the spectrum of $A:(a,b)\to \Cnn$. 
Then $\tau_0$ is a bifurcation point of the spectrum of $A$. In fact, in
every (open) neighborhood   $\EuScript{N}$ of $\tau_0$ there exist  bifurcation points $t_0$. Hence, the function $N$ is not constant on
$\EuScript{N}$. As this happens for every neighborhood of
$\tau_0$, we deduce that $\tau_0$ is a
bifurcation point.

Therefore, if we suppose that the bifurcation points of the spectrum of $A$ are isolated,
it will follow that the set of these points cannot have accumulation points that belong to $(a,b)$. But   
$a$ or $b$ may  be accumulation points of such set.

\smallskip

\begin{lemma} \label{lema2}
Let $A:(a,b)\to \Cnn$ be a continuous matrix function such that the bifurcation points of its spectrum, if any, are isolated.
 Let us denote  
\[
r:=\max_{t\in (a,b)}\#\, \s{A(t)}.       
\]
Then there are continuous functions
 $\lambda_1,\ldots,\lambda_r:(a,b)\to\C$ such
 that for all $t\in (a,b)$,     
\begin{equation}\label{eq:10}
\Lambda(A(t))=\{ \lambda_1(t),\ldots,\lambda_r(t)\}.
\end{equation}
If $t'_0 < t_0 < t''_0$, with $t'_0,t_0,t''_0 \in (a,b)$, are three consecutive  bifurcation points, then  
 $\#\, \Lambda(A(t))$ is constant in each  of the intervals
$(t'_0,t_0)$ and
 $(t_0,t''_0)$, let us  say that
\[
\#\, \Lambda(A(t))=\left\{ \begin{array}{l}
                        u \mbox{  if  } t\in (t'_0,t_0), \\
                        v \mbox{  if  } t\in (t_0,t''_0).
 \end{array}\right.
\]
\end{lemma}
\smallskip

\smallskip

   \emph{%
  Moreover, there exist two subsets  $\{j_1,\ldots,j_u\}$ and $\{k_1,\ldots,k_v\}$ of \underbar{$r$} such that for all   
 $t\in (t'_0,t_0)$ the numbers $\lambda_{j_1}(t),\ldots,\lambda_{j_u}(t)$
are the $u$  distinct eigenvalues of $A(t)$, and for all $t\in (t_0,t''_0)$ the numbers
$\lambda_{k_1}(t),\ldots,\lambda_{k_v}(t)$ are the $v$ distinct eigenvalues of $A(t)$.}

\smallskip

\demo
By Lemma~\ref{lem:maximo}, there are $r$ continuous functions  $\lambda_1,
\ldots,\lambda_r :(a,b)\to \C$ such that for all $t\in (a,b)$,
\[
\Lambda(A(t))=\{\lambda_1(t),\ldots,\lambda_r(t)\}.
\]
Given that $t'_0$ are $t_0$ two consecutive bifurcation points of the spectrum  of $A$, we deduce that $\Lambda(A(t))$ has a constant number of elements when $t$ runs over $(t'_0,t_0)$; this is consequence of Lemma~\ref{lem:local_constant};
let us say that $\#\, \Lambda(A(t))=u$ for  $t\in (t'_0,t_0)$. By Lemma~\ref{lem:nobifurca}  
if two of the functions $\lambda_1,\ldots,\lambda_r$, let us say $\lambda_i$ and
$\lambda_j$, coincide at a point  $\tau_0$ of $(t'_0,t_0)$
\[
\lambda_i(\tau_0)=\lambda_j(\tau_0),
\]
then $\lambda_i(t)=\lambda_j(t)$ for all $t\in (t'_0,t_0)$. This let us  select the
subset $\{\lambda_{j_1},\ldots,\lambda_{j_u}\}$ of $\{\lambda_1,\ldots,\lambda_r\}$ constituted by
the $u$ functions such that for all $t\in (t'_0,t_0)$,
\[
\Lambda(A(t))=\{\lambda_{j_1}(t),\ldots,\lambda_{j_u}(t)\}.
\]
The rest of the proof is analogous.\qed

The behavior  of the spectrum at a point that is not a bifurcation point, is more simple, as the following lemma shows. See a proof in Lemma 4.2, p. 163, of~\cite{evard:1985}.

\smallskip

\begin{lemma} \label{lem:90-segre-10}
Let $A:(a,b)\to \mathbb C^{n\times n}$ be a continuous matrix function. Let
$$r:=\max_{t\in(a,b)}\enspace \#\,\Lambda(A(t)).$$
If $t_0\in(a,b)$ is not a bifurcation point of the spectrum of $A$, then there are a number $\varepsilon >0$ and $q$ continuous functions  $\ffi_1,\ldots,\ffi_q:(a,b)\to\C$, $q\le r$, such that for  $t\in(t_0-\varepsilon,t_0+\varepsilon)$
\[
\Lambda \big(A(t)\big)=\{\ffi_1(t),\ldots,\ffi_q(t)\}    
\]
and the numbers $\ffi_1(t),\ldots,\ffi_q(t)$ are distinct.   
\end{lemma}

\subsection*{Constant multiplicities}

The following  lemma, proved in~\cite[Lemma 2.1, p. 365]{evard_gr:1990} with the help of Lemmas~\ref{lem:local_constant} and~\ref{lem:nobifurca}, says that if a monic polynomial has continuous coefficients
$c_i(t)$, and a constant number of distinct roots when $t$ varies in an interval of the real line, then the multiplicities of its  roots are constant also.

\smallskip
\begin{lemma} \label{lem:constancia}
Let $I$ be an interval of $\R$. Let $n,r$ be positive integers with $r\le n$. Let 
\[
c_1,\ldots,c_n \colon I\to\C
\]
be continuous functions.  Suppose that the number of distinct roots, $r$, of the polynomial
\[
p_t(\la)=\la^n+c_1(t)\la^{n-1}+\cdots +c_{n-1}(t)\la+c_n(t),\quad \la\in\C.
\]
is constant. Then there are $r$ continuous functions $\secsim \al r\colon I \to \C$
such that for all $t\in I$, the values
$\tresvari \al t r$ are the roots of $p_t(\la)=0$.
Moreover, the multiplicities of these roots are constant.
\end{lemma}

\smallskip

\section{Multiplicities of structured pseudoeigenvalues} 

By Theorem~\ref{th:uno} each connected component of the $\ep$-pseu\-do\-spec\-trum of structure $S$ of $A$ contains an eigenvalue of $A$; which implies that the number of connected components is $\le n$. The main theorem in this paper is the following result, which asserts that the sum of algebraic multiplicities of the $\ep$-pseudoeigenvalues of structure $S$ of $A$ associated with $z,\norma{z}<\ep,$ that are inside the same connected component, remains constant.

\smallskip

\begin{theorem} \label{th:conserva}
Let $T$ be a connected component of $\La'_{S,\ep}(A)$. Then for all $z\in\C^s$ such that $\norma{z}<\ep $,
\begin{equation}\label{eq:conserva}
\sum_{\zeta\in\s{A+M_S(z)}\cap T}
\m(\zeta, A+M_S(z))=\sum_{\al\in\s A \cap T}\m(
\al,A).
\end{equation} 
\end{theorem}

\smallskip

\demo  For all $t\in\R$ let us define
$Z(t):=A+M_S(tz)$. For all $t\in[0,1]$,
\[
\norma{tz}=\abso{t}\,\norma{z}< \ep;
\]
which implies $\s{Z(t)}\subset \La'_{S,\ep}(A)$. Even more, given that the ball $B(0,\ep)\subset \C^s$ is an open set, 
there  exists a $\de >0$ such that $\forall t\in (-\de,1+\de )$,
\[
\s{Z(t)}\subset \La'_{S,\ep}(A).
\]
It is obvious that $Z(0)=A$ y $Z(1)=A+M_S(z)$. 
The eigenvalues of $Z(t)$ describe continuous trajectories 
\[
\tresvari \la t n
\]
that go from
$\tresvari \la 0 n$ to $\tresvari \la 1 n$ inside $\La'_{S,\ep}(A)$, because $[0,1]$ is connected and the continuous image of a connected set is connected. We are going to see that in this trip of the eigenvalues, the sum of multiplicities of the ones lying inside each connected component of $ \La'_{S,\ep}(A)$ is preserved.

For each real $t$, let
\[
p_{Z(t)}(\la):=\det (\la I-Z(t))
\]
be the characteristic  polynomial of $Z(t)$; obviously,
\[
p_{Z(t)}(\la)=\polimonivar \la a n t,
\]
where each  of the functions
 $\tresvari a t n$ is a polynomial in $t$. Let
\[
u:=\max_{t\in \R}\#\s{Z(t)}.
\]
It is evident that $u$ may be less that $n$. In general, we have the next assertion.
 
 \textbf{Assertion 1.-} There exists a finite subset (or empty) $F$ of $\R$ such that
\[
\card \s{Z(t)}\begin{cases}
=u \text{ if } t\in \R\setminus F, \\
<u \text{ if } t\in F.
\end{cases}
\]
\textbf{Proof.} Let
\[
R(t):=R\big(p_{Z(t)},p'_{Z(t)}\big)
\]
be the Sylvester resultant matrix of the polynomials  
 $p_{Z(t)}$ and its derivative
 $p'_{Z(t)}$. The matrix $R(t)$, given by
\[
\!\!\!\!\!\!\begin{pmatrix}
 1 & a_1(t)   & \cdots & a_n(t) & 0 &\cdots & 0\\
 0 & 1 & a_1(t) & \cdots  & a_n(t) & \cdots & 0\\
 \vdots & & \ddots & \ddots & &\ddots &\vdots\\
 0 & 0 & \cdots & 1 & a_1(t)& \cdots   & a_n(t)\\
 n & (n-1)a_1(t)  &  \cdots & a_{n-1}(t) & 0 &\cdots & 0\\
 0 & n & (n-1)a_1(t)  &\cdots  & a_{n-1}(t) & \cdots & 0\\
 \vdots & & \ddots & \ddots & &\ddots &\vdots\\
 0 & 0 & \cdots & n & (n-1)a_1(t)& \cdots   & a_{n-1}(t)
   \end{pmatrix},
\]
is square of size $2n-1$.

Let $k, (n\le k\le 2n-1)$, be the rank of $R(t)$ as polynomial matrix in $t$; that is to say, $k$ is the greatest order of minors of $R(t)$ that are not identically null.
 Let us denote by $\tresvari D t \ell$ the minors of order $k$ that are not identically null.
Hence, if $\tau$ is a real value such that some of the numbers
$\tresvari D \tau \ell $ is different from 0, it follows $\rk R(\tau)=k$. Let $F$ be the set of the common real roots of the polynomials
\[
\tresvari D t \ell.
\]
Then, if $t_0\in F$, we have
$\rk R(t_0)<k$.
Therefore, if $t\in \R\setminus F$, then
$\nu \big(R(t)\big)=2n-1-k$. By Corollary~\ref{cor:distintas}, 
\[
u=n-\nu\big(R(t)\big)=n-(2n-1-k)=k+1-n.
\]
If $t_0\in F$, then $\nu \big(R(t)\big)>2n-1-k$;
whence
\[
n-\nu\big(R(t)\big)<n-(2n-1-k)=u.\qed
\]

By Definition~\ref{bifur}, $F$ is the set 
of bifurcation points of the spectrum of the matrix function
$Z\colon \R \to \Cnn$. There are only two significative cases in the relation of $F$ with the interval $[0,1]$:
\begin{itemize}
\item Case 1.- either, $F\cap [0,1]\neq \o$,
\item Case 2.- or, $F\cap [0,1]= \o$.
\end{itemize}

\textbf{Case 1.-} Let $t_1:=\min F\cap [0,1]$. Let $t_2$ be
the next element of $F$ that follows to $t_1$ (it can happen that
$t_2\in [0,1]$, or  $1<t_2$, or $F\cap [0,1]=\{t_1\}$).

  Suppose that $0\le t_1<t_2\le 1$. In virtue of Lemma~\ref{lem:nobifurca}, because $(t_1,t_2)$  is connected and
\[
\{\tresvari \la t n\}
\]
has exactly  $u$ elements for  $t\in (t_1,t_2)$,
there is a subset of functions
$
\{\secutres \la ju \}$  of  $ \{\secsim \la n\}$
such that for  $t\in (t_1,t_2)$,
\[
\{\cuatrovari \la jtu\}=\s{Z(t)}.
\]
As the functions $\secsim \la n$ are continuous,
\[
\s{Z(t_1)}=\{\cuatrovari \la j{t_1}u\};
\]
but, taking into account that $\s{Z(t_1)}$ has less than $u$ elements, there are repetitions in the list $\cuatrovari \la j{t_1}u$.

Let $T$ be a connected component of $\La'_{S,\ep}(A)$. Suppose that 
\[
\s{Z(t_1)}\cap T=\{\secsim \al g\};
\]
then there is a subset $\{\secsim{\tilde{\la}}p\}$ of $\{\secutres \la ju\}$ such that:
\begin{enumerate}
\item \begin{align*}
 \al_1 &= \tilde{\la}_1(t_1) =\cdots = \tilde{\la}_{\ga_1}(t_1),  \\
\al_2 &= \tilde{\la}_{\ga_1 +1}(t_1)= \cdots = \tilde{\la}_{\ga_1+\ga_2}(t_1),  \\
   \phantom{rrrrrrrrrr} &\phantom{sssssssssss}\vdots  \\
\al_g &= \tilde{\la}_{\ga_1 +\cdots +\ga_{g-1}+1}(t_1)= \\
\cdots & = \tilde{\la}_{\ga_1+\cdots+\ga_{g-1}+\ga_g}(t_1),  
\end{align*}
with $\ga_1+\cdots+\ga_g=p \quad(\Rightarrow g\le p)$;

\item 
\[
\{\secsim \al g\}\cap \{\tresvari{\la'}{t_1}{u-p}\}=\o,
\]
where
\[
\{\secsim{\la'}{u-p}\}:=\{\secutres \la ju\}\setminus \{\secsim{\tilde{\la}}p\}.
\]

\end{enumerate}

If $\de >0$ is sufficiently small, by Lemma~\ref{lem:cont_espectro}, for all $t\in (t_1,t_1+\de)$,
\begin{align*}
\m(\al_1,Z(t_1))&=\sum_{i=1}^{\ga_1}\m(\tilde{\la}_i
(t),Z(t)),\\
  &\phantom{s} \vdots\\
\m(\al_g,Z(t_1))&=\sum_{i=1}^{\ga_g}\m(\tilde{\la}_{\ga_1+\cdots +\ga_{g-1}+i}(t),Z(t));
\end{align*} 
hence, agreeing that $\sum_{i=1}^0 \ga_i:=0$,
\begin{multline}\label{eq:1}
\sum_{j=1}^g \m(\al_j,Z(t_1))=
\sum_{j=1}^g \sum_{k=1}^{\ga_j} \m(\tilde{\la}_{\sum_{i=1}^{j-1}\ga_i+k}(t),Z(t))\\
=\sum_{\ell=1}^p\m(\tla_\ell (t),Z(t)).
\end{multline}
By Lemma~\ref{lem:constancia} and due to $(t_1,t_2)$ is connected, the multiplicities
\begin{equation}\label{eq:t1t2}
\m(\la_{j_1}(t),Z(t)),\ldots,\m(\la_{j_u}(t),Z(t))
\end{equation}
are constant on $(t_1,t_2)$. 

Let
\[
\{\secsim \be h\}:=\s{Z(t_2)}\cap T;
\]
then by continuity we have
\[
\{\secsim \be h\}=\{\tresvari {\tla}{t_2}p\};
\]
in the set of the right hand side of this equality can have repetitions; so, let \{\secsim{\ttla}p\} be a reordering of $\{\secsim{\tla}p \}$ such that

\begin{align*}
 \be_1 &= \ttla_1(t_2) =\cdots = \ttla_{\eta_1}(t_2),  \\
\be_2 &= \ttla_{\eta_1 +1}(t_2)= \cdots = \ttla_{\eta_1+\eta_2}(t_2),  \\
   \phantom{rrrrrrrrrr} &\phantom{sssssssssss}\vdots  \\
\be_h &= \ttla_{\sum_{i=1}^{h-1}\eta_i+1}(t_2)= \cdots\\
& = \ttla_{\sum_{i=1}^h \eta_i}(t_2),  
\end{align*}
with $\eta_1+\cdots+\eta_h=p \quad(\Rightarrow h\le p)$. By~\eqref{eq:t1t2} and invoking again Lemma~\ref{lem:cont_espectro}, for  $t\in(t_1,t_2)$,
\begin{equation}\label{eq:2}
\sum_{j=1}^h \m(\be_j,Z(t_2))=\sum_{\ell=1}^p \m(\ttla_\ell (t),Z(t)).
\end{equation}
From~\eqref{eq:1} and~\eqref{eq:2},
\[
\sum_{j=1}^g \m(\al_j,Z(t_1))=\sum_{j=1}^h \m(\be_j,Z(t_2));
\]
which is equivalent to
\begin{equation}\label{eq:3}
\sum_{\al \in \s{Z(t_1)}\cap T}\m(\al,Z(t_1))
=\sum_{\be \in \s{Z(t_2)}\cap T}\m(\be,Z(t_2)).
\end{equation}

Now let $t_q$ be the last element of $F$ that is $\le 1$; then, iterating  $q-1$ times the preceding  reasoning, we arrive to 
\begin{equation}\label{eq:4}
\sum_{\al \in \s{Z(t_1)}\cap T}\m(\al,Z(t_1))
=\sum_{\be \in \s{Z(t_q)}\cap T}\m(\be,Z(t_q)).
\end{equation}
If it were true that $t_1=0$ and $t_q=1$, we would have finished the proof because  $Z(0)=A$ and $Z(1)=A+M_S(z)$. If it were true that $0<t_1$, let $t_0<0$ be the element of $F$ immediately  less than $t_1$ (if $t_0$ does  not exist, we take $t_0=-\infty $); then the number of distinct eigenvalues of $Z(t)$ when $t$ runs over $(t_0,t_1)$ would be constant and equal to $u$; in fact, they should be the elements of the set
\[
\{\cuatrovari \la itu\}
\]
for certain subset $\{\secutres \la iu\}$ of
$\{\secsim \la n\}$. 

Let 
\[
\{\secsim \hla f\}\subset \{\secutres \la iu\}
\]
be such that
\[
\s{A}\cap T=\{\tresvari {\hla}0f\}\quad (\text{it cannot have repetitions})
\]
with
\[
\s{A}\cap T\cap \{\tresvari {\la''}{0}{u-f}\}=\o
\]
and
\[
\{\secsim{\la''}{u-f}\}:=\{\secutres \la iu\}\setminus\{\secsim \hla f\}.
\]
Taking into account that the $u$ algebraic multiplicities
\[
\m(\la_{i_a}(t),Z(t)),\qquad \reco a1u
\]
are constant on $(t_0,t_1)$, for  $t$ in this interval,
\begin{multline*}
\m(\hla_1(t),Z(t))+\cdots+\m(\hla_f(t),Z(t))=
\sum_{i=1}^f \m(\hla_i(0),Z(0))=\\
\sum_{
\al\in\s{A}\cap T} \m(\al,A);
\end{multline*}
which implies for $ t\in(t_0,t_1)$,
\begin{multline*}
\sum_{
\al\in\s{A}\cap T} \m(\al,A)=\sum_{i=1}^f \m(\hla_i(t),Z(t))=
\sum_{j=1}^g \m(\al_j,Z(t_1));
\end{multline*}
this, and~\eqref{eq:4}, imply
\begin{equation}\label{eq:5}
\sum_{
\al\in\s{A}\cap T} \m(\al,A)=\sum_{
\be\in\s{Z(t_q)}\cap T} \m(\be,Z(t_q)).
\end{equation} 
If it were true that $t_q=1$, we would have finished the proof. If it were true that $t_q<1$, let $t_{q+1}$ be the element of $F$ that follows to $t_q$ (if it does not exist, we take $t_{q+1}=\infty$). Then, there is a subset
$\{\secutres \la ku\}$ of $\{\secsim \la n\}$ that parameterizes the $u$ eigenvalues of $Z(t)$ on the interval $(t_q,t_{q+1})$. 

Let 
\[
\{\secsim \hhla s\}\subset \{\secutres \la ku\}
\]
be such that
\[
\s{X}\cap T=\{\tresvari {\hhla}1s\}\quad (\text{it cannot have repetitions})
\]
with
\[
\s{X}\cap T\cap \{\tresvari {\la'''}{1}{u-s}\}=\o
\]
and
\[
\{\secsim{\la'''}{u-s}\}:=\{\secutres \la ku\}\setminus\{\secsim \hhla s\}.
\]
Then, as the values
\[
\m(\hhla_i(t),Z(t))\qquad \reco i1s
\]
are constant on $(t_q,t_{q+1})$, 
\begin{equation}\label{eq:6}
\sum_{
\be\in\s{Z(t_q)}\cap T} \m(\be,Z(t_q))=\sum_{
\xi\in\s{X)}\cap T} \m(\xi,X).
\end{equation}
Due to~\eqref{eq:5} and~\eqref{eq:6}, we have proved that
\begin{equation}\label{eq:7}
\sum_{
\al\in\s{A)}\cap T} \m(\al,A)=\sum_{
\xi\in\s{X)}\cap T} \m(\xi,X).
\end{equation}

\textbf{Case 2.-} If $F\cap [0,1]=\o$, the equality~\eqref{eq:7} follows easily. If $\{\secutres \la ju\}$
 is the subset of $\{\secsim \la n\}$ that parameterizes the eigenvalues of $Z(t)$ from $t=0$ to $t=1$, we deduce that each of the multiplicities $\m(\la_{j_i}(t),Z(t)),\; \reco i1u$, is constant, and in consequence for $ \reco i1u$,
\[
\m(\la_{j_i}(0),Z(0))=\m(\la_{j_i}(1),Z(1));
\]
this proofs~\eqref{eq:7} in this case.\qed

\section{Consequences}
\begin{remark}\label{obs:raizcomun}
{\em If $t_0\in F$, then
\[
D_1(t_0)=0,\ldots, D_\ell(t_0)=0.
\]
Therefore, the maximum multiplicity of the roots of $p_{Z(t_0)}(\la)$ is greater than or equal to the maximum multiplicity of the roots of $p_{Z(t)}(\la)$ for $t\notin F$;
\[
\max_{\al\in \s {Z(t)}} \m(\al, Z(t))\le\max_{\mu\in \s {Z(t_0)}} \m(\mu, Z(t_0)).
\]
But the number
\[
\max_{\mu\in \s {Z(t_0)}} \m(\mu, Z(t_0))
\]
can change really when  $t_0$ varies in $F$; this occurs if
 all the minors  not identically  null of $R(t)$ of orders $k-1, k-2, \ldots,$ vanish   at $t_0$.
}
\end{remark}

\begin{remark}\label{obs:factores}
{\em Let $d_j(t)$ be the $j$th determinantal divisor  of the polynomial matrix $R(t)$; that is,  $d_j(t)$
is the greatest common divisor of all minors of order
 $j$ of $R(t), \reco j1k$; in particular,
\[
d_k(t)=\gcd \{\tresvari D t \ell\}
\]
If $f_1(t)|\cdots | f_k(t)$ are the nonzero  invariant factors of  $R(t)$, we have for each $\reco j1k$,
\[
f_j(t)=\frac{d_j(t)}{d_{j-1}(t)}, \text{ with } d_0(t):=1.
\]
Hence
\[
f_1(t)\cdots  f_j(t)=d_j(t);
\]
in particular $f_1(t)\cdots  f_k(t)=d_k(t)$. Whence, 
we conclude that $F$ is the set of real roots of $d_k(t)=0$,
\[
F=\{t_0\in \R \mid d_k(t_0)=0 \}.
\]
}
\end{remark}
\begin{remark}\label{obs:duda}
{\em Briefly, let us denote $p_t(\la):=p_{Z(t)}(\la)$ and by
\[
R_q(t):=R(p_t,p'_t,\ldots,p^{(q)}_t)
\]
 the generalized Sylvester resultant matrix of the polynomials in $\la$: $p_t,p'_t,\ldots,p^{(q)}_t$,
$\reco q1n$.
Let $k_q=\rk R_q(t)$ be as matrix of polynomials in $t$.
With this notation, the preceding resultant matrix $R(t)$ is $R_1(t)$ and $k_1=k$. Let $d^{(q)}_{k_q}(t)$  be the $k_q$th determinantal divisor   of $R_q(t)$. Denote by $F_q$ the set of real roots of $d^{(q)}_{k_q}(t)=0$.
What relations exist among the sets $F_q $ and $F$?
}
\end{remark}

With more detail, let $\varrho_q(t)$ be the number of roots of multiplicity $q$ of the characteristic polynomial $p_{Z(t)}(\la)$, for $\reco q1n$. 
With the notations of the proof of Theorem~\ref{th:conserva},
for  $t\in \R\setminus F$,
\[
u=\varrho_1(t)+\cdots+\varrho_n(t)
\]

We will say that $t_0\in \R$ is a bifurcation point of multiplicities of this polynomial, if the function
\[
t\mapsto \big(\tresvari {\varrho}tn \big)
\]
is not constant on every neighborhood of $t_0$. Let us denote by $G$ the set of these points. In the proof of Theorem~\ref{th:conserva}
we have demonstrated that the points $t$ of $\R \setminus F$ are not bifurcation points of multiplicities. Hence, $G\subset F$.
May not be empty the set $G$?

\begin{corollary}\label{cor:nearest}
Let $A\in \Cnn$. If $\ep>0$ is such  that $\La'_{S,\ep}(A)$ has $n$ connected components, then for all vector $z$ that satisfies $\norma{z}<\ep$, the matrix $A+M_S(z)$ has simple  eigenvalues.  
\end{corollary}

For all matrix $X\in \Cnn$, let us denote by
\[
p_X(\la)
\]
its characteristic polynomial. Denote by
$R(p_X,p'_X,\ldots,p^{(k)}_X)$ the generalized Sylvester resultant matrix of $p_X$ and its successive derivatives
$p_X',\ldots,p^{(k)}_X$.
The polynomial admits a root of multiplicity at least $k+1$ if and only if the nullity of the resultant matrix is greater than 0
\[
\nu \big(R(p_X,p'_X,\ldots,p^{(k)}_X)\big)>0.
\]

\smallskip

\begin{corollary}\label{cor:nearest2}
Let $A\in \Cnn$. If $\ep>0$ is such that $\La'_\ep(A)$  has $\rho(\ep)$ connected components $\secsim S{\rho(\ep)}$. Let
\[
\mu(\ep,A):=\max_{1\le i\le \rho(\ep)} \sum_{\al\in \s A \cap 
S_i} \m(\al,A).
\] 
Then every matrix $X$ such that $\norma{X-A}<\ep$ has its eigenvalues of multiplicity  $\le\mu(\ep,A)$.
\end{corollary}

\subsection{A lower bound}

\begin{corollary}\label{cor:ultimo}
Let $\ep>0$ and let $
\secsim S {\rho(\ep)}$ be the connected components of $\La'_\ep(A)$. Let us define
\begin{align*}
\mu(\ep,A)&:=\max_{1\le i\le \rho(\ep)} \sum_{\al\in \s A \cap 
S_i} \m(\al,A),\\
\m(X)&:=\max_{\xi\in\s X} \m(\xi,X), \text{ where } X\in\Cnn.
\end{align*}
Let $k$  be an integer such that $\m(A)\le k\le n.$
Then
\begin{equation*}
\sup \;\{\ep>0\colon \mu(\ep,A)\le k\}
\le \min_{\m(X)\ge k+1}\|X-A\|.
\end{equation*}
\end{corollary}

\subsection*{ Acknowledgments} 

I thank  Ion Zaballa for the idea of the proof of Theorem~\ref{th:uno}.

\end{document}